\documentclass[12pt]{article}

\usepackage{amsmath}
\usepackage{amssymb}
\usepackage{amsfonts}
\usepackage{amsopn}
\usepackage[all]{xy}
\usepackage{amsthm}
\usepackage{amscd}

\swapnumbers
\theoremstyle{plain}
\newtheorem{thm}{Theorem}[section]
\newtheorem{lem}[thm]{Lemma}
\newtheorem{cor}[thm]{Corollary}

\theoremstyle{definition}
\newtheorem{ntt}[thm]{}
\newtheorem{rem}[thm]{Remark}

\newcommand{\pl}{\mathbb{P}}   
\newcommand{\zz}{\mathbb{Z}}   
\newcommand{\F}{{\mathrm{F}_4}}
\newcommand{\G}{\mathrm{G}_2}  
\newcommand{\B}{\mathrm{B}}
\newcommand{\w}{\bar{\omega}} 
\newcommand{\id}{\mathrm{id}}       
\newcommand{\pr}{\mathit{pr}}       

\newcommand{\Var}{\mathcal{V}ar} 
\newcommand{\M}{\mathcal{M}}     
\newcommand{\Cor}{\mathcal{C}or} 
\newcommand{\zab}{\mathbb{Z}\text{-}\mathcal{A}b} 

\newcommand{\hra}{\hookrightarrow}     
\newcommand{\lra}{\longrightarrow}     
\newcommand{\xra}[1]{\xrightarrow{#1}} 

\DeclareMathOperator{\Pic}{\mathrm{Pic}}    
\DeclareMathOperator{\Br}{\mathrm{Br}}     
\DeclareMathOperator{\Gal}{\mathrm{Gal}}   
\DeclareMathOperator{\CH}{\mathrm{CH}}      
\DeclareMathOperator{\Mor}{\mathrm{Mor}}
\DeclareMathOperator{\End}{\mathrm{End}}
\DeclareMathOperator{\im}{\mathrm{Im}}      
\DeclareMathOperator{\Ker}{\mathrm{Ker}}      
\DeclareMathOperator{\K}{\mathrm{K}}        
\DeclareMathOperator{\HD}{\mathcal{H}}      

\title{Motivic decomposition of anisotropic varieties of
type $\F$ into generalized Rost motives
}

\author{S.~Nikolenko, N.~Semenov, K.~Zainoulline  \footnote{
Supported partially by DAAD, INTAS 00-566, Alexander von Humboldt foundation,
RTN-Network HPRN-CT-2002-00287.}}
\date{}

\begin{document}

\maketitle

\begin{abstract}
We prove that the Chow motive of 
an anisotropic projective homogeneous variety of type $\F$
is isomorphic to the direct sum of twisted copies 
of a generalized Rost motive. 
In particular, we provide an explicit construction of a generalized Rost
motive for a generically splitting variety for a symbol in $\K^M_3(k)/3$.
We also establish a motivic isomorphism between 
two anisotropic non-isomorphic projective
homogeneous varieties of type $\F$.
All our results hold for Chow motives with integral coefficients.
\end{abstract}

\section{Introduction}

The subject of the present paper begins
with the celebrated result of M.~Rost \cite{Ro98}
devoted to the motivic decomposition of a norm quadric.  
The existence of such a decomposition became one of the main
ingredients in the proof of Milnor conjecture by V.~Voevodsky. 
The generalization of this conjecture to other primes $p>2$,
known as the Bloch-Kato conjecture, was proven recently by M.~Rost and
V.~Voevodsky. One of the ingredients of the proof is the fact
that the motive with $\zz/p\zz$-coefficients of a splitting (norm) variety $X$ contains as a direct
summand a certain geometric motive $M_{p-1}$ 
called \emph{generalized Rost motive} \cite[Sect.~5]{Vo03}.
This motive is indecomposable and 
splits as a direct sum of twisted Lefschetz motives over the separable
closure of the base field.

Note that Voevodsky's construction of $M_{p-1}$ 
relies heavily on the language of triangulated
category of motives.
The main goal of the present paper
is to provide an explicit and shortened construction
of the motive $M_2$ ($p=3$)  
working only within the classical category of Chow motives. 
More precisely, we provide such a construction 
for an exceptional projective homogeneous variety of type $\F$ which
splits the symbol (in $\K^M_3(k)/3$) given by the Rost-Serre invariant $g_3$ 
(see~\ref{normvar}). 
 
Note that if $X$ is generically cellular and splits a pure symbol,
it is expected that the motive of $X$
is isomorphic to the direct sum of twisted copies of the motive $M_{p-1}$.
The motivic decomposition that we obtain confirms these expectations.
Namely, we prove the following

\begin{thm}\label{mainthm} 
Let $X$ be an anisotropic variety over a 
field $k$ such that over a cubic field extension $k'$ of $k$ 
it becomes isomorphic to the projective homogeneous variety $G/P$,
where $G$ is a split simple group of type $\F$ and 
$P$ its maximal parabolic subgroup corresponding
to the last or the first three vertices of the Dynkin diagram 
\ref{DynkinF4}.

Then the Chow motive of $X$ (with integral coefficients) is isomorphic to
the direct sum of twisted copies of an indecomposable motive $R$
\begin{equation}\label{maindec}
\M(X)\cong \bigoplus_{i=0}^7 R(i)
\end{equation}
which has the property that
over $k'$ it
becomes isomorphic to the direct sum of twisted Lefschetz motives 
$\zz\oplus \zz(4)\oplus \zz(8)$.
\end{thm}   

The next result provides the first
known ``purely exceptional''  example of two
different anisotropic varieties with isomorphic motives.
Recall that the similar result
for groups of type $\mathrm{G}_2$ obtained in \cite{Bo03} 
provides a motivic isomorphism
between quadric and an exceptional Fano variety.

\begin{thm}\label{motisom}
Let $\mathcal{G}$ be an anisotropic algebraic group over $k$
such that over a cubic field extension $k'$ of $k$ it becomes
isomorphic to a split simple group $G$ of type $\F$.
Let $X_1$ and $X_4$ be two projective
$\mathcal{G}$-homogeneous varieties which over $k'$
become isomorphic to $G/P_1$
and $G/P_4$ respectively, where $P_1$ (resp. $P_4$) is the maximal parabolic subgroup corresponding to the last (resp. first)
three vertices of the Dynkin diagram \ref{DynkinF4}.

Then the motives of $X_1$ and $X_4$ are isomorphic.
\end{thm}

The main motivation for our work was the result of N.~Karpenko where
he gave a shortened construction of a Rost motive for a norm quadric 
\cite{Ka98}. The key idea is to produce enough idempotents 
in the ring $\CH(X\times X)$ considered over the separable closure of $k$
and then lift them to $k$ using the
Rost Nilpotence Theorem (see \cite{CGM}).
Contrary to the techniques used by Voevodsky, 
the proof of Theorem \ref{mainthm} is based on well-known and elementary facts
about linear algebraic groups, projective homogeneous varieties,
and Chow groups. 

We expect that similar methods can be applied
to other projective homogeneous varieties,
thus providing analogous motivic decompositions.
In particular, applying our arguments to a Pfister quadric
one obtains the celebrated decomposition 
into Rost motives (see \cite[Example~7.3]{KM02}).
For exceptional groups of type $\G$ one immediately obtains
the motivic decomposition of the Fano variety
together with the motivic isomorphism constructed in \cite{Bo03}.

The paper is organized as follows. 
In Section \ref{ratcyc} we provide background
information on Chow motives and rational cycles. 
Section \ref{hadiag} is devoted to
computational matters of Chow rings. Namely, we introduce Pieri and 
Giambelli formulae and discuss their relationships with Hasse diagrams.
In Section \ref{multtables} we apply the formulae introduced in 
Section \ref{hadiag}
to projective homogeneous varieties $X_1$ and $X_4$ of type $\F$. 
In Section \ref{constproj} we prove Theorems \ref{mainthm} and \ref{motisom}.
Finally, in the Appendix we explain the intermediate technical steps
of computations used in the proofs.

\section{Motives and rational cycles}\label{ratcyc}

In the present section 
we introduce the category of Chow motives over a field $k$ 
following \cite{Ma68} and \cite{Ka01}. We remind the notion of 
a rational cycle and state the Rost Nilpotence
Theorem for idempotents following \cite{CGM}.

\begin{ntt} Let $k$ be a field and $\Var_k$ be a category
of smooth projective varieties over $k$. 
First, we define the category of \emph{correspondences} (over $k$)
denoted by $\Cor_k$. 
Its objects are smooth projective varieties over $k$. 
For morphisms, called correspondences,
we set $\Mor(X,Y):=\CH^{\dim X}(X\times Y)$.
For any two correspondences $\alpha\in \CH(X\times Y)$ and 
$\beta\in \CH(Y\times Z)$ we define their composition 
$\beta\circ\alpha\in \CH(X\times Z)$ as
\begin{equation}\label{compos}
\beta\circ\alpha ={\pr_{13}}_*(\pr_{12}^*(\alpha)\cdot \pr_{23}^*(\beta)),
\end{equation}
where $\pr_{ij}$ denotes the projection
on the $i$-th and $j$-th factors of $X\times Y\times Z$ 
respectively and ${\pr_{ij}}_*$, ${\pr_{ij}^*}$ denote
the induced push-forwards and pull-backs for Chow groups. 

The pseudo-abelian completion of $\Cor_k$ is called the category
of \emph{Chow motives} and is denoted by $\M_k$.
The objects of $\M_k$
are pairs $(X,p)$, where $X$ is a smooth projective variety
and $p$ is an idempotent, that is, $p\circ p=p$. 
The morphisms between two objects $(X,p)$ and $(Y,q)$ 
are the compositions 
$q\circ \Mor(X,Y) \circ p$.
The motive $(X,\id)$ will be denoted by $\M(X)$.
\end{ntt}

\begin{ntt}
By construction, $\M_k$ is a tensor additive category with self-duality,
where the self-duality is given by the transposition of cycles $\alpha\mapsto \alpha^t$,
and the tensor product is given by the usual fiber product
$(X,p)\otimes(Y,q)=(X\times Y, p\times q)$.
Moreover, the contravariant Chow functor 
$\CH: \Var_k \to \zab$ (to the category of $\zz$-graded abelian groups)
factors through $\M_k$, that is, one has a commutative diagram of functors
$$
\xymatrix{
\Var_k \ar[rr]^{\CH}\ar[rd]_{\Gamma}& & \zab \\
 & \M_k \ar[ru]_{CH}& 
}
$$
where $\Gamma: f\mapsto \Gamma_f$ is the covariant graph functor and 
$CH:(X,p) \mapsto \im(p^*)$ is the realization.
\end{ntt}

\begin{ntt} Observe that the composition product $\circ$
induces the ring structure on the abelian group $\CH^{\dim X}(X\times X)$.
The unit element of this ring is the class of the diagonal map $\Delta_X$, which is defined by
$\Delta_X\circ \alpha = \alpha\circ \Delta_X=\alpha$ for all 
$\alpha\in \CH^{\dim X}(X\times X)$.
\end{ntt}

\begin{ntt} Consider the morphism
$(e,\id):\{pt\} \times \pl^1 \to \pl^1\times\pl^1$. 
Its image by means of
the induced push-forward $(e,\id)_*$ does not depend on the choice of the point
$e:\{pt\}\to \pl^1$
and defines the projector in $\CH^1(\pl^1\times\pl^1)$ denoted by $p_1$.
The motive $\zz(1)=(\pl^1,p_1)$ is called \emph{Lefschetz motive}.
For a motive $M$ and a nonnegative integer $i$ 
we denote its twist by $M(i)=M\otimes \zz(1)^{\otimes i}$.
\end{ntt}

\begin{ntt}\label{compcy} 
Let $G$ be a split simple linear algebraic group over $k$.
Let $X$ be a projective $G$-homogeneous variety, that is, 
$X\simeq G/P$, where $P$ is a parabolic subgroup of $G$.
The abelian group structure of $\CH(X)$, as well as
its ring structure, is well-known.
Namely, $X$ has a cellular filtration and 
the generators of Chow groups of
the bases of this filtration correspond to
the free additive generators of $\CH(X)$.
Note that the product of two projective homogeneous varieties
$X\times Y$ has a cellular filtration as well,
and $\CH^*(X\times Y)\cong \CH^*(X)\otimes \CH^*(Y)$
as graded rings.
The correspondence product of two cycles 
$\alpha=f_\alpha \times g_\alpha \in \CH(X\times Y)$ and
$\beta=f_\beta \times g_\beta \in \CH(Y\times X)$ is given
by (cf. \cite[Lem.~5]{Bo03})
\begin{equation}\label{composecycles}
(f_\beta\times g_\beta)\circ(f_\alpha\times g_\alpha)=
\deg(g_\alpha \cdot f_\beta)(f_\alpha\times g_\beta),
\end{equation}
where $\deg: \CH(Y)\to \CH(\{pt\})=\zz$ is the degree map.
\end{ntt}

\begin{ntt}
Let $X$ be a projective variety of dimension $n$ over a field $k$. 
Let $k'$ be a field extension of $k$ and $X'=X\times_k k'$.
We say a cycle $J\in \CH(X')$ is {\it rational}
if it lies in the image of the pull-back homomorphism
$\CH(X)\to \CH(X')$.
For instance, there is an obvious rational cycle $\Delta_{X'}$ on 
$\CH^n(X'\times X')$ that is given by the diagonal class.
Clearly, all linear combinations, 
intersections and correspondence products of rational cycles
are rational.
\end{ntt}

\begin{ntt}
Several techniques allow
to produce rational cycles 
(cf. \cite[Prop.~3.3]{Ka04} for the case of quadrics).
We shall use the following:
\begin{itemize}
\item[(i)]
Consider a variety $Y$ and a morphism $X\to Y$ such that
$X'=Y'\times_Y X$, where $Y'=Y\times_k k'$. Then any rational cycle
on $\CH(Y')$ gives rise to a rational cycle on $\CH(X')$ by 
the induced pull-back $\CH(Y')\to \CH(X')$.

\item[(ii)]
Consider a variety $Y$ and a projective morphism
$Y\to X$ such that $Y'=X'\times_X Y$. Then any rational cycle
on $\CH(Y')$ gives rise to a rational cycle on $\CH(X')$
by the induced push-forward $\CH(Y')\to \CH(X')$.

\item[(iii)]
\label{genpoint}
Let $X$ and $Y$ be projective homogeneous varieties over $k$
such that $X$ splits completely (i.e., the respective group
splits) over the function field $k(Y)$.
Consider the following pull-back diagram
$$
\xymatrix{
\CH^i(X\times Y)\ar[d]_{f} \ar[r]^g& \CH^i(X'\times Y')\ar[d]^{f'}\\
\CH^i(X_{k(Y)})\ar[r]^{=} & \CH^i(X'_{k'(Y')})
}
$$
where the vertical arrows are surjective by \cite[\S5]{IK00}.
Now take any cycle $\alpha\in \CH^i(X'\times Y')$, $i\leq \dim X$.
Let $\beta=g(f^{-1}(f'(\alpha)))$.
Then $f'(\beta)=f'(\alpha)$ and $\beta$ is rational.
Hence, $\beta=\alpha+J$, where $J\in \Ker f'$, and we conclude
that $\alpha+J\in \CH^i(X'\times Y')$ is rational.
\end{itemize}
\end{ntt}

\begin{ntt}[Rost Nilpotence]\label{exproj} 
Finally, we shall also use the following fact (see \cite[Cor.~8.3]{CGM})
that follows from Rost Nilpotence Theorem. 
Let $p'$ be a non-trivial rational idempotent on $\CH^n(X'\times X')$,
i.e., $p'\circ p'=p'$. Then there exists a non-trivial idempotent $p$
on $\CH^n(X\times X)$ such that $p\times_k k'=p'$. 
Hence, the existence of a non-trivial rational idempotent $p'$ on 
$\CH^n(X'\times X')$ gives rise to
the decomposition of the Chow motive of $X$
$$
\M(X)\cong(X,p)\oplus (X,\id_X-p).
$$
\end{ntt}



\section{Hasse diagrams and Chow rings}\label{hadiag}

To each projective homogeneous variety $X$ we may associate
an oriented labeled graph $\HD$ called Hasse diagram. 
It is known that the ring structure of $\CH(X)$ is determined by $\HD$.
In the present section we remind several facts concerning
relations between Hasse diagrams and Chow rings. For detailed explanations
of these relations see \cite{De74}, \cite{Hi82a} and \cite{Ko91}.

\begin{ntt}
Let $G$ be a split simple algebraic group defined over a field $k$.
We fix a maximal split torus $T$ in $G$ and a Borel subgroup $B$ of $G$
containing $T$ and defined over $k$. We denote by $\Phi$ the root system of $G$,
by $\Pi$ the set of simple roots of $\Phi$ corresponding to $B$, by $W$ the Weyl group, and by $S$
the corresponding set of fundamental reflections.

Let $P=P_\Theta$ be a (standard) parabolic subgroup corresponding
to a subset $\Theta\subset\Pi$, i.e., $P=BW_\Theta B$, where
$W_\Theta=\langle s_\theta, \theta\in\Theta\rangle$.
Denote $$W^\Theta=\{w\in W\mid\forall\, s\in\Theta\quad l(ws)=l(w)+1\},$$
where $l$ is the length function.
The pairing 
$$W^\Theta\times W_\Theta\to W\qquad (w,v)\mapsto wv$$
is a bijection and $l(wv)=l(w)+l(v)$.
It is easy to see that $W^\Theta$ consists of all representatives
in the cosets $W/W_\Theta$ which have minimal length.
Sometimes it is also convenient to
consider the set of all representatives of maximal length. We shall
denote this set as ${}^\Theta W$. Observe that there is a 
bijection $W^\Theta \to {}^\Theta W$ given by 
$v\mapsto vw_\theta$, where $w_\theta$ is the longest element of $W_\Theta$.
The longest element of $W^\Theta$ corresponds to the
longest element $w_0$ of the Weyl group.
\end{ntt}

\begin{ntt}\label{defHasse}
To a subset $\Theta$ of the finite set $\Pi$ 
we associate an oriented labeled graph,
which we call a Hasse diagram and denote by $\mathcal{H}_W(\Theta)$.
This graph is constructed as follows.
The vertices of this graph are the elements of $W^\Theta$. 
There is an edge from a vertex $w$ to a vertex $w'$
labelled with $i$ if and only if $l(w)<l(w')$ and $w'=ws_i$.
A sample Hasse diagram is provided in \ref{Hasse}.
Observe that the diagram $\mathcal{H}_W(\emptyset)$ 
coincides with the Cayley graph associated to the pair $(W,S)$.
\end{ntt}

\begin{lem}
The assignment $\mathcal{H}_W\colon\Theta \mapsto \mathcal{H}_W(\Theta)$
is a contravariant functor from the category of subsets 
of the finite set $\Pi$
(with embeddings as morphisms) to the category of oriented graphs. 
\end{lem}

\begin{proof}
It is enough to embed the diagram $\mathcal{H}_W(\Theta)$ 
to the diagram $\mathcal{H}_W(\emptyset)$.
We do this as follows.
We identify the vertices of $\mathcal{H}_W(\Theta)$
with the subset of vertices of $\mathcal{H}_W(\emptyset)$
by means of the bijection $W^\Theta\to {}^\Theta W$.
Then the edge
from $w$ to $w'$ of ${}^\Theta W\subset W$ has a label $i$ if and only if
$l(w)<l(w')$ and $w'=s_iw$ (as elements of $W$). 
Clearly, the obtained graph will coincide
with $\mathcal{H}_W(\Theta)$.
\end{proof}

\begin{ntt}\label{ChowHasse} Now consider the Chow ring of a projective
homogeneous variety $G/P_\Theta$. 
It is well known that $\CH(G/P_\Theta)$ is a free abelian group
with a basis given by varieties $[X_w]$ that correspond 
to the
vertices $w$ of the Hasse diagram $\mathcal{H}_W(\Theta)$. 
The degree of the basis element $[X_w]$ 
corresponds to the
minimal number of edges needed to connect the respective vertex $w$
with $w_\theta$ (which is the longest word).
The multiplicative structure of $\CH(G/P_\Theta)$ depends
only on the root system of $G$ and the diagram $\mathcal{H}_W(\Theta)$. 
\end{ntt}

By definition one immediately obtains

\begin{lem} The contravariant functor 
$\CH\colon \Theta \mapsto \CH(G/P_\Theta)$
factors through the category of Hasse diagrams $\mathcal{H}_W$, i.e.,
the pull-back (ring inclusion)
$$\CH(G/P_{\Theta'})\hra \CH(G/P_\Theta)$$ arising from 
the embedding $\Theta \subset \Theta'$ is induced by
the embedding of the respective Hasse diagrams 
$\mathcal{H}_W(\Theta') \subset \mathcal{H}_W(\Theta)$.
\end{lem}

\begin{cor}\label{embed} Let $B$ be a Borel subgroup of $G$
and $P$ its (standard) parabolic subgroup.
Then $\CH(G/P)$ is a subring of $\CH(G/B)$. 
The generators of $\CH(G/P)$ are $[X_w]$, where
$w\in{^\Theta W}\subset W$. 
The cycle  $[X_w]$ in
$\CH(G/P)$ has the codimension $l(w_0)-l(w)$.
\end{cor}

\begin{proof} Apply the lemma to the case $B=P_\emptyset$ and $P=P_{\Theta'}$.
\end{proof}

Hence, in order to compute $\CH(G/P)$
it is enough to compute $\CH(X)$, where $X=G/B$ is the variety of
complete flags.
The following results provide tools to perform such computations.

\begin{ntt}\label{complem}
In order to multiply two basis elements $h$ and $g$ of $\CH(G/P)$ such that 
$\deg h + \deg g = \dim G/P$ we use the following formula 
(see \cite[1.4]{Ko91}):
$$
[X_w]\cdot [X_{w'}]= \delta_{w,w_0w'w_\theta}\cdot [pt].
$$
\end{ntt}

\begin{ntt}[Pieri formula]\label{Pieri}
In order to multiply two basis elements of $\CH(X)$ one of which
is of codimension 1 we use the following formula 
(see \cite[Cor.~2 of 4.4]{De74}):
$$
[X_{w_0s_\alpha}][X_w]=
\sum_{\beta\in\Phi^+,\, l(ws_\beta)=l(w)-1}
\langle\beta^\vee,\overline{\omega}_\alpha\rangle[X_{ws_\beta}],
$$
where the sum runs through
the set of positive roots $\beta\in \Phi^+$, $s_\alpha$ denotes
the simple reflection corresponding to $\alpha$
and $\w_\alpha$ is the fundamental weight corresponding
to $\alpha$. Here $[X_{w_0s_\alpha}]$ is the element of codimension 1.
\end{ntt}

\begin{ntt}[Giambelli formula]\label{Giambelli}
Let $P=P(\Phi)$ be the weight space. 
We denote as $\w_1,\ldots\w_l$ the basis of $P$
consisting of fundamental weights. 
The symmetric algebra $S^*(P)$ is isomorphic to
$\mathbb{Z}[\w_1,\ldots\w_l]$. 
The Weyl group $W$ acts on $P$, hence, on $S^*(P)$. 
Namely, for a simple root $\alpha_i$,
$$w_{\alpha_i}(\w_j)=
\begin{cases}
\w_i-\alpha_i, & i=j, \\
\w_j, & \text{otherwise}.
\end{cases}$$
We define a linear map $c\colon S^*(P)\to\CH^*(G/B)$ as follows.
For a homogeneous $u\in\mathbb{Z}[\w_1,\ldots,\w_l]$ 
$$
c(u)=\sum_{w\in W,\, l(w)=\deg(u)}\Delta_w(u)[X_{w_0w}],
$$ 
where for $w=w_{\alpha_1}\ldots w_{\alpha_k}$ we denote by $\Delta_w$
the composition of derivations
$\Delta_{\alpha_1}\circ\ldots\circ\Delta_{\alpha_k}$ and
the derivation $\Delta_{\alpha_i}\colon S^*(P)\to S^{*-1}(P)$ is 
defined by $\Delta_{\alpha_i}(u)=\frac{u-w_{\alpha_i}(u)}{\alpha_i}$.
Then (see \cite[ch.~IV, 2.4]{Hi82a})
$$
[X_w]=c(\Delta_{w^{-1}}(\frac{d}{|W|})),
$$ 
where $d$ is the product
of all positive roots in $S^*(P)$. In other words, the element 
$\Delta_{w^{-1}}(\frac{d}{|W|})\in c^{-1}([X_w])$.

Hence, in order to multiply two basis elements $h,g\in\CH(X)$
one may take their preimages under the map $c$
and multiply them in $S^*(P)=\mathbb{Z}[\w_1,\ldots\w_l]$, 
finally applying $c$ to the product.
\end{ntt}



\section{Homogeneous varieties of type $\F$}\label{multtables}

In the present section we remind several well-known facts  
concerning Albert algebras, groups of type $\F$ and 
respective projective homogeneous varieties
(see  \cite{PR94} and \cite{Inv}).
At the end we provide partial computations of Chow rings of these varieties.

\begin{ntt}\label{DynkinF4}
From now on 
let $X_i$ be an anisotropic projective homogeneous variety over a field $k$ 
that over a cubic field extension $k'$ of $k$ 
becomes isomorphic
to the projective homogeneous variety $X'_i=G/P_i$, where
$G$ is a split group of type $\F$ and
$P_i=P_{\theta_i}$ 
is a maximal parabolic subgroup of $G$
corresponding to the subset $\theta_i=\{1,2,3,4\}\setminus \{i\}$
of the Dynkin diagram
$$
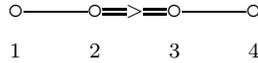
\xymatrix@M=0pt{
\circ \ar@{}[d] |-{1} \ar@{-}[r]& \circ\ar@{}[d] |-{2}\ar@{=}[r] |-{>}& \circ\ar@{}[d] |-{3} \ar@{-}[r]& \circ\ar@{}[d] |-{4} \\
 & & &
}
$$
\end{ntt}

\begin{rem}
The variety $X_i$
is a $_\xi{G}$-variety over $k$, where 
$_\xi{G}$
is an algebraic group over $k$ that is the
twisted form of a split group $G$ of type $\F$ by means of a $1$-cocycle
$\xi\in H^1(k,G(k'))$ (see \cite[Prop.~4]{De77}). 

If the base field $k$ has characteristic not $3$ and $\mu_3\subset k$,
then all such groups $_\xi{G}$ are automorphism 
groups of Albert algebras coming from the first Tits construction.
\end{rem}

The next two important properties of the varieties $X_i$ will be extensively
used in the sequel.

\begin{lem}\label{ratPicard}
The Picard group $\Pic(X'_i)$
is a free abelian group of rank 1
with a rational generator.
\end{lem}

\begin{proof}
Since $P$ is maximal, $\Pic(X'_i)$
is a free abelian group of rank 1.
Consider the following exact sequence (see \cite{Ar82} and \cite[2.3]{MT95})
$$ 
0 \lra \Pic X_i \lra (\Pic X'_i)^\Gamma \xra{\alpha} \Br(k),
$$
where
$\Gamma=\Gal(k'/k)$ is the Galois group
and $\Br(k)$ the Brauer group of $k$. 
The map $\alpha$ is explicitly 
described in \cite{MT95} in terms of Tits classes.
Since all groups of type $\F$ are simply-connected and adjoint
their Tits classes are trivial and so is $\alpha$. 
Since $\Gamma$ acts trivially on $\Pic(X'_i)$ and 
the image of $\alpha$ is trivial,
we have $\Pic(X_i)\simeq \Pic(X'_i)$.
\end{proof}

\begin{lem}\label{albertsplit}
For any $i,j\in \{1,2,3,4\}$, the variety $X_i$ splits completely
over the function field $k(X_j)$.
\end{lem}

\begin{proof} The following arguments belong to S.~Garibaldi.

It is equivalent to show that the function field $K=k(X_j)$ 
splits the group $\widetilde{G}={}_\xi{G}$.
First, observe that if $\widetilde{G}$ is isotropic then it is split.
Indeed, by Tits classification \cite[p.~60]{Ti66}, the only
other possibility is that $\widetilde{G}$ has a maximal parabolic $P_4$ 
defined over $k$, but no others. That is, the semisimple anisotropic
kernel $H$ of $G$ is of type ${\B}_3$. Since $\widetilde{G}$ is split
by a cubic field extension, the same is true for $H$. But this
is impossible for an anisotropic group of type ${\B}_3$
by Springer's Theorem on quadratic forms and odd-degree extensions.

Now let $k'$ be a cubic extension of $k$ that splits $\widetilde{G}$.
Since $K$ is a regular extension of $k$, the tensor product
$K'=K \otimes_k k'$ is a field and has dimension 1 or 3 over $K$.
Hence, $\widetilde{G}\otimes_k K$ is isotropic and 
is split by the extension $K'$. By previous arguments
$\widetilde{G}\otimes_k K$ is split.
\end{proof}

\begin{rem}\label{normvar} 
To any Albert algebra $A$
one can associate a symbol $\{a_1,a_2,a_3\}$ in $\K^M_3(k)/3$
given by Rost-Serre's invariant $g_3$.
If $A$ is obtained by the first Tits construction,
then $g_3=0$ iff $A$ is split,
and the variety $X_j$ corresponding to $A$ provides an example of a 
generically splitting variety
for the symbol given by $g_3$ (see \cite[Def.~1.8]{Su05}).
Observe that $X_j$ is not exactly a norm ($\nu_n$-)variety 
in the sense of \cite[Sect.~4]{Vo03}, 
since it has the wrong dimension,
but rather a ``Pfister quadric'' version of it. 
\end{rem}

\begin{ntt}\label{XY}
From now on we consider only two anisotropic
varieties $X_1$ and $X_4$. Recall
that over a cubic field extension $k'$ they become isomorphic to
$X_1'=G/P_1$ and $X_4'=G/P_4$ 
respectively (see \ref{DynkinF4}).
Varieties $X_1'$ and $X_4'$ are not isomorphic 
(their Chow rings are different) and 
have the dimension $15$. Observe that 
the variety $X_4$ is a twisted
form of a hyperplane section of the Cayley plane
$\mathbb{OP}^2$ which was extensively studied in \cite{IM05}.
\end{ntt}

\begin{ntt}\label{Hasse} 
The Hasse diagram (see \ref{defHasse}) for $X_1'$  looks as follows
$$
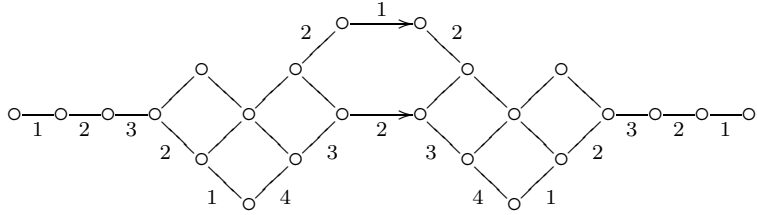
\xymatrix@M=0pt@=1em{
 & & & & & & &\circ\ar[rr]^1 & & \circ \ar@{-}[rd]^2& & & & & & & \\
 & & & & \circ \ar@{-}[rd]& & \circ \ar@{-}[rd]\ar@{-}[ru]^2& & & & \circ \ar@{-}[rd]& & \circ \ar@{-}[rd] & & & & \\
\circ \ar@{-}[r]_1 & \circ \ar@{-}[r]_2 & \circ \ar@{-}[r]_3 & \circ \ar@{-}[ru]\ar@{-}[rd]_2& & \circ \ar@{-}[ru]\ar@{-}[rd]& & \circ  \ar[rr]_2&  & \circ \ar@{-}[ru]\ar@{-}[rd]_3& & \circ \ar@{-}[ru]\ar@{-}[rd]& & \circ \ar@{-}[r]_3 & \circ \ar@{-}[r]_2 & \circ \ar@{-}[r]_1 & \circ \\
 & & & & \circ \ar@{-}[ru]\ar@{-}[rd]_1 & & \circ\ar@{-}[ru]_3 & & & & \circ\ar@{-}[ru]\ar@{-}[rd]_4 & & \circ \ar@{-}[ru]_2& & & & \\
 & & & & & \circ \ar@{-}[ru]_4& & &  & & & \circ\ar@{-}[ru]_1 & & & & & 
}
$$
and for $X_4'$
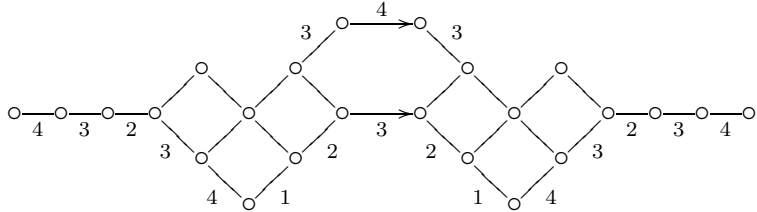
$$
\xymatrix@M=0pt@=1em{
 & & & & & & &\circ\ar[rr]^4 & & \circ \ar@{-}[rd]^3& & & & & & & \\
 & & & & \circ \ar@{-}[rd]& & \circ \ar@{-}[rd]\ar@{-}[ru]^3& & & & \circ \ar@{-}[rd]& & \circ \ar@{-}[rd] & & & & \\
\circ \ar@{-}[r]_4 & \circ \ar@{-}[r]_3 & \circ \ar@{-}[r]_2 & \circ \ar@{-}[ru]\ar@{-}[rd]_3& & \circ \ar@{-}[ru]\ar@{-}[rd]& & \circ  \ar[rr]_3&  & \circ \ar@{-}[ru]\ar@{-}[rd]_2& & \circ \ar@{-}[ru]\ar@{-}[rd]& & \circ \ar@{-}[r]_2 & \circ \ar@{-}[r]_3 & \circ \ar@{-}[r]_4 & \circ \\
 & & & & \circ \ar@{-}[ru]\ar@{-}[rd]_4 & & \circ\ar@{-}[ru]_2 & & & & \circ\ar@{-}[ru]\ar@{-}[rd]_1 & & \circ \ar@{-}[ru]_3& & & & \\
 & & & & & \circ \ar@{-}[ru]_1& & &  & & & \circ\ar@{-}[ru]_4 & & & & & 
}
$$
We draw the diagrams in such a way that 
the labels on opposite sides of a parallelogram are equal, and in that
case we omit all labels but one.
Recall that (see \ref{ChowHasse}) the vertices of this graph correspond to the basis elements of the Chow group.
The leftmost vertex is the unit class and the rightmost one 
is the class of a $0$-cycle of degree $1$. 
\end{ntt}

\begin{ntt}\label{notbas}
We denote the basis elements of the respective
Chow groups as follows
$$
\CH^i(X_1')=
\begin{cases}
\langle h_1^i\rangle,& i=0,\ldots,3,12,\ldots,15,\\
\langle h_1^i, h_2^i \rangle,& i=4,\ldots,11.  
\end{cases}$$
$$
\CH^i(X_4')=
\begin{cases}
\langle g_1^i\rangle,& i=0,\ldots,3,12,\ldots,15,\\
\langle g_1^i, g_2^i \rangle,& i=4,\ldots,11.  
\end{cases}$$

The generators with lower index $1$ correspond 
to the lower vertices of the respective Hasse diagrams,
and with index $2$ to the upper ones (if there are two generators).
\end{ntt}

\begin{ntt}\label{comdim}
Applying \ref{complem} we immediately obtain the following partial 
multiplication table:
$$
h_i^s h_j^{15-s} = \delta_{ij}h_1^{15}, 
\qquad g_i^s g_j^{15-s} = \delta_{ij} g_1^{15},
$$
where $\delta_{ij}=1$ if $i=j$ and $0$ otherwise.
\end{ntt}

\begin{ntt}\label{Pieritab}
By Pieri formula \ref{Pieri} 
we obtain the following partial multiplication tables
for $\CH(X_1')$:
         \begin{alignat*}{5}
         h_1^1h_1^1 &=h_1^2, &\quad h_1^1h_1^2 &=2h_1^3,
         &\quad h_1^1h_1^3&=2h_2^4+h_1^4, &
         \quad h_1h_2^4 & =h_2^5, \\ h_1^1h_1^4 &=2h_2^5+h_1^5,&\quad
         h_1^1h_2^5 & =2h_2^6+h_1^6, &\quad h_1^1h_1^5 & =2h_1^6,
         &\quad h_1^1h_2^6&=h_2^7+h_1^7, \\
         h_1^1h_1^6 & =2h_1^7, &\quad h_1^1h_2^7 &=2h_2^8+h_1^8, &\quad
         h_1^1h_1^7 & =h_2^8+2h_1^8, &\quad h_1^1h_2^8 & =h_2^9, \\
         h_1^1h_1^8&=h_2^9+2h_1^9, & \quad h_1^1h_2^9 & =2h_2^{10},
         &\quad h_1^1h_1^9 &=h_2^{10}+2h_1^{10}, &\quad
         h_1^1h_2^{10} & =h_2^{11}+2h_1^{11}, \\
         h_1^1h_1^{10} & =h_1^{11}, &\quad h_1^1h_2^{11}&=2h_1^{12}, &
         \quad h_1^1h_1^{11} & =h_1^{12}, &\quad h_1^1h_1^{12} &=2h_1^{13},\\
         h_1^1h_1^{13} & =h_1^{14}, &\quad h_1^1h_1^{14} & =h_1^{15}.
         \end{alignat*}
for $\CH(X_4')$:
\begin{alignat*}{5}
g_1^1g_1^1 &=g_1^2, &\quad g_1^1g_1^2 &=g_1^3,
&\quad g_1^1g_1^3&=g_1^4+g_2^4, &
\quad g_1^1g_2^4 & =g_2^5, \\ g_1^1g_1^4 &=g_1^5+g_2^5,&\quad
g_1^1g_2^5 & =g_1^6+g_2^6, &\quad g_1^1g_1^5 & =g_1^6,
&\quad g_1^1g_2^6&=g_1^7+g_2^7, \\
g_1^1g_1^6 & =g_1^7, &\quad g_1^1g_1^7 &=2g_1^8+g_2^8, &\quad
g_1^1g_2^7 & =g_1^8+2g_2^8, &\quad g_1^1g_2^8 & =g_2^9, \\
g_1^1g_1^8&=g_1^9+g_2^9, & \quad g_1^1g_2^9 & =g_2^{10},
&\quad g_1^1g_1^9 &=g_1^{10}+g_2^{10}, &\quad
g_1^1g_2^{10} & =g_1^{11}+g_2^{11}, \\
g_1^1g_1^{10} & =g_1^{11}, &\quad g_1^1g_2^{11}&=g_1^{12}, &
\quad g_1^1g_1^{11} & =g_1^{12}, &\quad g_1^1g_1^{12} &=g_1^{13},\\
g_1^1g_1^{13} & =g_1^{14}, &\quad g_1^1g_1^{14} & =g_1^{15}.
\end{alignat*}
\end{ntt}

\begin{ntt}\label{graphpieri}
Observe that the multiplication tables \ref{Pieritab} 
can be visualized by means of
slightly modified Hasse diagrams.
Namely, for the variety $X_1'$ consider the following graph which
is obtained from the respective Hasse diagram by adding a few more edges
and erasing all the labels:
$$
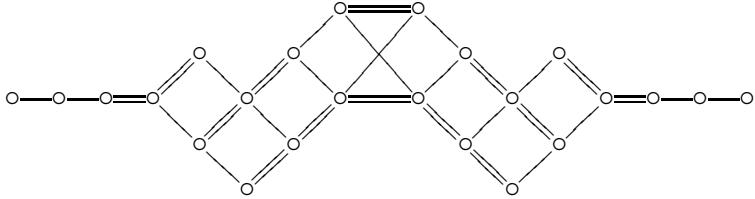
\xymatrix@M=0pt@=1em{
 & & & & & & &\circ\ar@{=}[rr]\ar@{-}[rrdd] & & \circ \ar@{-}[rd]\ar@{-}[lldd] & & & & & & & \\
 & & & & \circ \ar@{-}[rd]& & \circ \ar@{-}[rd]\ar@{-}[ru]& & & & \circ \ar@{=}[rd]& & \circ \ar@{=}[rd] & & & & \\
\circ \ar@{-}[r] & \circ \ar@{-}[r] & \circ \ar@{=}[r] & \circ \ar@{=}[ru]\ar@{-}[rd]& & \circ \ar@{=}[ru]\ar@{-}[rd]& & \circ  \ar@{=}[rr]&  & \circ \ar@{-}[ru]\ar@{=}[rd]& & \circ \ar@{-}[ru]\ar@{=}[rd]& & \circ \ar@{=}[r] & \circ \ar@{-}[r] & \circ \ar@{-}[r] & \circ \\
 & & & & \circ \ar@{=}[ru]\ar@{-}[rd] & & \circ\ar@{=}[ru] & & & & \circ\ar@{-}[ru]\ar@{=}[rd] & & \circ \ar@{-}[ru]& & & & \\
 & & & & & \circ \ar@{=}[ru]& & &  & & & \circ\ar@{-}[ru] & & & & & 
}
$$
and for $X_4'$:
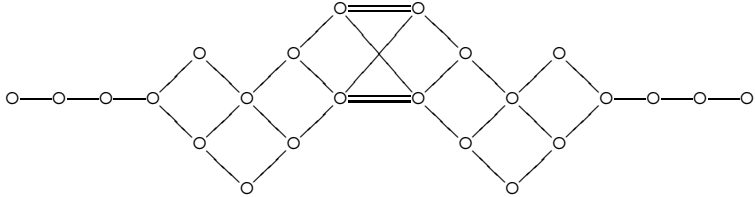
$$
\xymatrix@M=0pt@=1em{
 & & & & & & &\circ\ar@{=}[rr]\ar@{-}[rrdd] & & \circ \ar@{-}[rd]\ar@{-}[lldd] & & & & & & & \\
 & & & & \circ \ar@{-}[rd]& & \circ \ar@{-}[rd]\ar@{-}[ru]& & & & \circ \ar@{-}[rd]& & \circ \ar@{-}[rd] & & & & \\
\circ \ar@{-}[r] & \circ \ar@{-}[r] & \circ \ar@{-}[r] & \circ \ar@{-}[ru]\ar@{-}[rd]& & \circ \ar@{-}[ru]\ar@{-}[rd]& & \circ  \ar@{=}[rr]&  & \circ \ar@{-}[ru]\ar@{-}[rd]& & \circ \ar@{-}[ru]\ar@{-}[rd]& & \circ \ar@{-}[r] & \circ \ar@{-}[r] & \circ \ar@{-}[r] & \circ \\
 & & & & \circ \ar@{-}[ru]\ar@{-}[rd] & & \circ\ar@{-}[ru] & & & & \circ\ar@{-}[ru]\ar@{-}[rd] & & \circ \ar@{-}[ru]& & & & \\
 & & & & & \circ \ar@{-}[ru]& & &  & & & \circ\ar@{-}[ru] & & & & & 
}
$$
Then the multiplication rules 
can be restored as follows: 
for a vertex $u$ 
(that corresponds to a basis element of the Chow group) 
we set 
$$
H\cdot u = \sum_{u\to v} v,
$$ 
where $H$ denotes either $h_1^1$ or $g_1^1$ and 
the sum runs through all the edges going from $u$ 
one step to the right (cf. \cite[Cor.~3.3]{Hi82b}). 
\end{ntt}

\begin{ntt}\label{gprod}
Applying Giambelli formula \ref{Giambelli}
we obtain the following products (for details see Appendix):
$$
h_1^4h_1^4=8h_1^8+6h_2^8,\qquad g_1^4g_1^4=4g_1^8+3g_2^8.
$$ 
\end{ntt}


\section{Construction of rational idempotents}\label{constproj}

The goal of the present section is to prove Theorems \ref{mainthm} and
\ref{motisom}.
The proof consists of several steps. First, 
using properties \ref{ratPicard} and 
\ref{albertsplit} we provide several important cycles 
$\rho_i$ and prove their rationality. 
Multiplying and composing them, we obtain a set of
pairwise orthogonal idempotents $p'_i$ and $q'_i$.
Then, using Rost Nilpotence
Theorem (see \ref{exproj}) we obtain the desired motivic decomposition
and, hence, finish the proof of \ref{mainthm}.
At the end we construct
an explicit cycle which provides
a motivic isomorphism of Theorem \ref{motisom}.

\paragraph{\it Proof of Theorem \ref{mainthm}}

Let $X_1$ and $X_4$ be the varieties corresponding to the last (resp. first)
three roots of the Dynkin diagram (see \ref{XY}).
By the hypothesis of Theorem \ref{mainthm} the variety $X$ is isomorphic
either to $X_1$ or $X_4$ over $k$.
As in \ref{XY} let $k'$ denote the cubic field extension and
let $X_1'$ and $X_4'$ be the respective base change.
We start with the following obvious observation.

\begin{ntt}\label{sectrat} 
Since the variety $X$ of Theorem \ref{mainthm} splits by a cubic field
extension, 
transfer arguments show that any cycle of the kind 
$3z \in \CH(X')$ is rational.
Hence, to prove that a cycle in $\CH(X')$ is rational
it is enough to prove this modulo $3$.
We shall write $x\equiv y$ if $x-y=3z$ for some cycle $z$. 
\end{ntt}

\begin{ntt}\label{Pierirat} The rational cycles to start with 
one obtains by Lemma \ref{ratPicard}. Namely, 
those are the classes of rational generators of the Picard groups
$h_1^1$ and $g_1^1$ (see \ref{notbas}). 
Clearly, their powers
$(h_1^1)^i$ and $(g_1^1)^i$, $i=2,\ldots,7$ are rational as well.
\end{ntt}

\begin{ntt}\label{genrat} 
Apply the arguments of \ref{genpoint}(iii) 
to $\CH^4(X_1'\times X_4')$ 
(this can be done because of Lemma \ref{albertsplit}).
There exists a rational cycle $\alpha_1\in\CH^4(X_1'\times X_4')$ 
such that $f'(\alpha_1)=h_1^4\times 1$.
This cycle must have the following form:
$$ 
\alpha_1 = h_1^4\times 1 + a_1 h_1^3\times g_1^1 + 
a_2 h_1^2\times g_1^2 + 
a_3 h_1^1\times g_1^3 + a_4 1\times g_2^4 + a' 1\times g_1^4,
$$
where $a_i,a'\in\{-1,0,1\}$. 
We may reduce $\alpha_1$ by adding cycles that are known to be rational
by \ref{Pierirat} to 
$$
\alpha_1=(h_1^4\times 1) + a(1\times g_1^4),
$$
where $a\in\{-1,0,1\}$.
Repeating the same procedure for a rational cycle 
$\alpha_2 \in\CH^4(X_1'\times X_4')$
such that $f'(\alpha_2)=1\times g_1^4$ and reducing it we obtain
the rational cycle
$$
\alpha_2=b (h_1^4\times 1)+ (1\times g_1^4),
$$
where $b\in\{-1,0,1\}$.
Hence, there is a rational cycle of the form
$$
r = h_1^4\times 1 + \varepsilon \cdot (1\times g_1^4),
$$
where the (indefinite) coefficient $\varepsilon$ is either $1$ or $-1$.
\end{ntt}

Now combining \ref{Pierirat} and \ref{genrat} together we obtain

\begin{lem}\label{ratlem} For all $i=0,\ldots,7$ the cycles
$$
\rho_i=r^2\cdot ((h_1^1)^i\times (g_1^1)^{7-i}) \in \CH^{15}(X_1'\times X_4')
$$
are rational.
\end{lem}

\begin{ntt} \label{explproj} By the previous lemma all cycles
$$
\rho_{7-i}^t \circ \rho_i \in \CH^{15}(X_1'\times X_1') \text{ and }
\rho_i \circ \rho_{7-i}^t \in \CH^{15}(X_4'\times X_4')
$$
where $i=0,\ldots,3$ are rational.
Direct computations (see Appendix) show that these cycles are congruent
modulo 3 to the following cycles in $\CH^{15}(X_1'\times X_1')$:
\begin{align*}
\rho_7^t\circ \rho_0 &\equiv 
1\times h_1^{15}+h_1^4\times (h_1^{11}+h_2^{11})+ h_1^8\times (h_1^7+h_2^7), \\
\rho_6^t\circ\rho_1 &\equiv h_1^1\times h_1^{14}+
(2h_2^5+h_1^5)\times (h_2^{10}-h_1^{10})+
(h_2^9-h_1^9)\times h_2^6, \\
\rho_5^t\circ \rho_2 &\equiv h_1^2\times h_1^{13} + (h_1^6+h_2^6)\times h_1^9
+ (2h_1^{10}-h_2^{10})\times (h_1^5+h_2^5), \\
\rho_4^t\circ \rho_3 &\equiv  h_1^3\times h_1^{12} +
h_2^7\times (h_2^8-h_1^8) + h_2^{11}\times (h_2^4-h_1^4)
\end{align*}
and in $\CH^{15}(X_4'\times X_4')$:
\begin{align*}
\rho_0\circ \rho_7^t &\equiv 1\times g_1^{15}+ g_1^4\times (g_1^{11}-g_2^{11})+ g_1^8\times (g_1^7+g_2^7), \\
\rho_1\circ\rho_6^t &\equiv g_1^1\times g_1^{14}+(2g_1^5-g_2^5)\times (g_1^{10}+g_2^{10})+ (g_1^9+g_2^9)\times g_2^6, \\
\rho_2\circ\rho_5^t &\equiv g_1^2\times g_1^{13}+ (g_1^6-g_2^6)\times g_1^9 + (g_1^{10}+2g_2^{10})\times (g_2^5-g_1^5), \\
\rho_3\circ\rho_4^t &\equiv g_1^3\times g_1^{12} + g_2^7\times (g_2^8-g_1^8)+ g_2^{11}\times (g_1^4+g_2^4)
\end{align*}
respectively, which turn to be idempotents. 
Denote them
by $p_0'$, $p_1'$, $p_2'$, $p_3'$ and $q_0'$, $q_1'$, $q_2'$, $q_4'$ 
respectively. Since they are congruent to rational cycles,
they are rational.
To complete the picture observe that 
the transposed cycles $(p_i')^t$ and $(q_i')^t$, $i=0,\ldots,3$,
are rational idempotents as well.
Hence, we have produced eight rational idempotents in each Chow ring.
\end{ntt}

\begin{ntt}
Now easy computations show that these idempotents
are orthogonal to each other, and their sum is equal to the
respective diagonal cycle
\begin{align*}
\Delta_{X_1'} &=\sum_{i,j,s} \delta_{ij} h_i^s\times h_j^{15-s}=
\sum_{l=0}^3 p_l' + (p_l')^t \in \CH^{15}(X_1'\times X_1'), \\
\Delta_{X_4'} &=\sum_{i,j,s} \delta_{ij} g_i^s\times g_j^{15-s}=
\sum_{l=0}^3 q_l' + (q_l')^t \in \CH^{15}(X_4'\times X_4').
\end{align*}
Hence, by \ref{exproj} we obtain decompositions 
of the motives of $X_1$ and $X_4$
$$
\M(X_1)=\bigoplus_{i=0}^3 (X_1,p_i)\oplus(X_1,p_i^t),
$$
$$
\M(X_4)=\bigoplus_{i=0}^3 (X_4,q_i)\oplus(X_4,q_i^t),
$$
where $p_i\times_k k'=p_i'$ and $q_i\times_k k'=q_i'$
are pairwise orthogonal idempotents (over $k$).
\end{ntt}

\begin{ntt} By a straightforward computation using the definition
of the idempotents $p_i'$ and $q_i'$ given in \ref{explproj}
we immediately obtain that
\begin{itemize}
\item there are isomorphisms of motives
\begin{align*}
(X_1,p_i) &\simeq (X_1,p_0)(i), &\quad
(X_1,p_i^t) &\simeq (X_1,p_i)(7-i), \\
(X_4,q_i) &\simeq (X_4,q_0)(i), &\quad
(X_4,q_i^t) &\simeq (X_4,q_i)(7-i),
\end{align*}
\item over the cubic extension $k'$
the motives $(X_1,p_0)$ and $(X_4,q_0)$ 
split as direct sums of Lefschetz motives 
$\zz\oplus \zz(4)\oplus\zz(8)$, where 
the shifts correspond to the codimensions
of the first factors of $p_0$ and $q_0$.
\end{itemize}
\end{ntt}

\begin{ntt} To finish the proof of Theorem \ref{mainthm} we have to prove
that the motives $(X_1,p_0)$ and $(X_4,q_0)$ are indecomposable.
To see this observe that the group of endomorphisms $\End(X_1',p_0')$
is a free abelian group with the basis 
\begin{equation}\label{bas}
\langle 1\times h_1^{15}, 
h_1^4\times (h_1^{11}+h_2^{11}),h_1^8\times (h_1^7+h_2^7) \rangle
\end{equation}
Assume that $(X_1,p_0)$ is decomposable, 
then the motive $(X_1',p_0')$ is decomposable as well.
The latter means that there exists a non-trivial rational idempotent 
in $\End(X_1',p_0')$. 
In this case easy computations show that
one of the elements of the basis (\ref{bas}) must be rational.
For instance, assume $1\times h_1^{15}$ is rational. Then
the cycle $(1\times h_1^{15}) \cdot (h_1^{15}\times 1)=h_1^{15}\times h_1^{15}$ 
is rational and so is its image $h_1^{15}$ by means of the push-forward 
$\CH_0(X_1'\times X_1') \to \CH_0(X_1')$ induced by a projection.
Hence, we obtain a cycle of degree $1$ on $X_1$, i.e., the variety $X_1$
must have a rational point \cite[Cor.~p.~205]{PR94}. We have arrived to a contradiction. \qquad \qed
\end{ntt}

\begin{rem} The following picture demonstrates how the realizations
of motives $(X_1',p_i')$, $i=0,1,2,3$,
are supported by the generators of the Hasse diagram
(the numbers $i$ drawn inside the rectangulars correspond to the motives $(X_1',p_i')$).
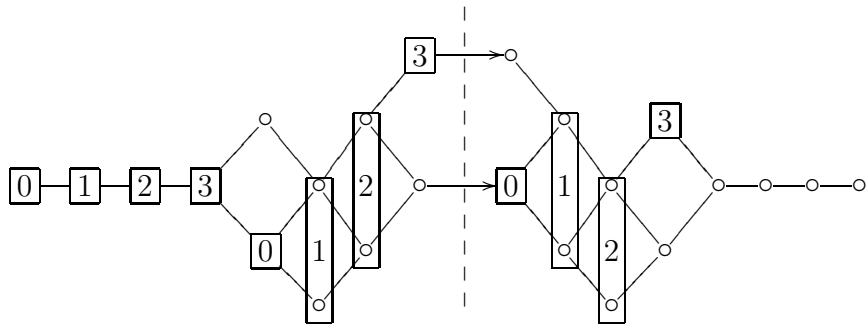
$$
\xymatrix@M=0pt@=1em{
 & & & & & & &                     & \ar@{--}[ddddd]               &                  & & & & & & & \\
 & & & & & & &*+<1ex>[F]{3}\ar[rr] &  & \circ \ar@{-}[rd]& & & & & & & \\
 & & & & \circ \ar@{-}[rd]& & \circ \ar@{-}[rd]\ar@{-}[ru]& & & & \circ \ar@{-}[rd]& & *+<1ex>[F]{3} \ar@{-}[rd] & & & & \\
*+<1ex>[F]{\small 0} \ar@{-}[r] & *+<1ex>[F]{1} \ar@{-}[r] & *+<1ex>[F]{2} \ar@{-}[r] & *+<1ex>[F]{3} \ar@{-}[ru]\ar@{-}[rd]& & \circ \ar@{-}[ru]\ar@{-}[rd]& 2 & \circ  \ar[rr]&  & *+<1ex>[F]{0} \ar@{-}[ru]\ar@{-}[rd]& 1 & \circ \ar@{-}[ru]\ar@{-}[rd]& & \circ \ar@{-}[r] & \circ \ar@{-}[r] & \circ \ar@{-}[r] & \circ \\
  &  &  &  & *+<1ex>[F]{0} \ar@{-}[ru]\ar@{-}[rd] & 1 & \circ\ar@{-}[ru] & & & & \circ\ar@{-}[ru]\ar@{-}[rd] & 2 & \circ \ar@{-}[ru]& & & & \\
 & & & & & \circ \ar@{-}[ru]& & &  & & & \circ\ar@{-}[ru] & & & & & 
\save "4,6"."6,6"*+<0.7ex>\frm{-}\restore \save "3,11"."5,11"*+<0.7ex>\frm{-}\restore
\save "3,7"."5,7"*+<0.7ex>\frm{-}\restore \save "4,12"."6,12"*+<0.7ex>\frm{-}\restore
}
$$
A similar picture for the motives $(X_1',p_i^t)$ corresponding to transposed idempotents
is obtained by the reflection along the vertical dashed line of the diagram.
\end{rem}

\begin{proof}[Proof of Theorem \ref{motisom}]
We use the notation of the proof of \ref{mainthm}.
Easy computations (see Appendix) show  that
$$
(\rho_0+\rho_1+\varepsilon \rho_2+\varepsilon \rho_3)+
(\rho_0+\rho_1+\varepsilon \rho_2+\varepsilon \rho_3)^t\equiv 
\sum_{i,j,s} \pm \delta_{ij} h_i^s \times g_j^{15-s}.
$$
Denote the right hand side by $J$.
We have
$J^t\circ J=\Delta_{X_1'}$ and $J\circ J^t=\Delta_{X_4'}$,
i.e., $J$ and $J^t$ are two mutually inverse rational correspondences and
the cycle $J$ provides
a rational motivic isomorphism between the motives of $X_1'$ and $X_4'$.  
By Rost Nilpotence Theorem, $J$
can be lifted to a motivic isomorphism between the motives of 
twisted forms $X_1$ and $X_4$.
This completes the proof of Theorem \ref{motisom}.
\end{proof}


\section*{Appendix}
Most of the computations of the present section
were performed and checked using the Maple package by J.~Stembridge \cite{St04}.

\begin{ntt} In the present paragraph 
we list the intermediate results of the computations
of \ref{explproj}. First, using \ref{gprod} we obtain
$$
r^2\equiv -h_1^8\times 1 -\varepsilon h_1^4\times g_1^4 + 1\times g_1^8.
$$
Then, using \ref{Pieritab} we obtain the following congruences
for the cycles $\rho_i$
\begin{align*}
\rho_0 &\equiv -1\times g_1^{15}+ \varepsilon h_1^4\times(g_1^{11}-g_2^{11})+ h_1^8\times (g_1^7+g_2^7) \\
\rho_1 &\equiv -h_1^1\times g_1^{14}+ \varepsilon (h_2^5-h_1^5)\times (g_1^{10}+g_2^{10}) +(h_2^9-h_1^9)\times g_2^6 \\
\rho_2 &\equiv -h_1^2\times g_1^{13} - \varepsilon (h_1^6+h_2^6)\times g_1^9 + (h_1^{10}+h_2^{10})\times (g_2^5-g_1^5) \\
\rho_3 &\equiv h_1^3\times g_1^{12} + \varepsilon h_2^7\times (g_2^8-g_1^8) - h_2^{11}\times (g_1^4+g_2^4) \\
\rho_4 &\equiv (h_1^4-h_2^4)\times g_2^{11} + \varepsilon (h_2^8-h_1^8)\times g_2^7 + h_1^{12}\times g_1^3 \\
\rho_5 &\equiv (h_1^5+h_2^5)\times (g_2^{10}-g_1^{10}) + \varepsilon h_1^9\times (g_2^6-g_1^6) - h_1^{13}\times g_1^2 \\
\rho_6 &\equiv h_2^6\times (g_1^9+g_2^9) + \varepsilon (h_2^{10}-h_1^{10})\times (g_1^5+g_2^5) - h_1^{14}\times g_1^1 \\
\rho_7 &\equiv (h_1^7+h_2^7)\times g_1^8 + \varepsilon (h_1^{11}+h_2^{11})\times g_1^4 - h_1^{15}\times 1
\end{align*}
\end{ntt}

\begin{ntt}
The root enumeration follows Bourbaki \cite{Bou}.
To obtain the products of \ref{gprod}
we apply the Giambelli formula \ref{Giambelli}.
Let $\w_i$, $i=1,\ldots,4$, be the fundamental weights.
Then, the preimages of $h_1^4$ and $g_1^4$
in $S^*(P)$ can be expressed as polynomials in fundamental weights 
as follows
\begin{multline*}
h_1^4=c(\frac{11}6 \w_1^2 \w_4^2+\frac34 \w_1^2\w_2^2 -\frac43 \w_1 \w_2 \w_3^2+
\frac{11}6 \w_1^2 \w_3^2-\frac23 \w_1 \w_2 \w_3 \w_4+\frac{11}{12} \w_1^4+\\
\frac16 \w_2^4-\frac43 \w_2\w_3^2 \w_4 +
\frac43 \w_2\w_3 \w_4^2 +\frac23 \w_2^2\w_3 \w_4 +\frac23  \w_1 \w_2\w_4^2-
\frac{11}6 \w_1^2\w_3 \w_4 +\\
2  \w_1 \w_3^2\w_4 - 2 \w_1 \w_3\w_4^2 -
\frac7{12} \w_1^3\w_2  -\frac{11}6 \w_1^2\w_2  \w_3+\frac43 \w_1 \w_2^2 \w_3+
\frac23 \w_2^2 \w_3^2 -\\
\frac23 \w_2^3 \w_3 -\frac13 \w_1\w_2^3  -\frac23 \w_2^2\w_4^2 ),
\end{multline*}
\begin{multline*}
g_1^4 = c(\frac{11}6\w_4^4-\frac76\w_3\w_4^3+\frac{11}{12}\w_1^2\w_4^2+
\frac32\w_3^2\w_4^2-\frac{11}6\w_2\w_3\w_4^2+\frac{11}{12}\w_2^2\w_4^2-\\
\frac{11}{12}\w_1\w_2\w_4^2-\frac23\w_3^3\w_4-\frac12\w_1^2\w_2\w_4+
\frac13\w_1^2\w_3\w_4+\frac43\w_2\w_3^2\w_4+\frac12\w_1\w_2^2\w_4-\\
\frac23\w_2^2\w_3\w_4-\frac13\w_1\w_2\w_3\w_4+\frac13\w_3^4-
\frac13\w_1\w_2^2\w_3+\frac13\w_1\w_2\w_3^2-\frac13\w_1^2\w_3^2+\\
\frac13\w_1^2\w_2\w_3+\frac13\w_2^2\w_3^2-\frac23\w_2\w_3^3).
\end{multline*}
Multiplying the respective polynomials and taking the $c$ function, 
we find the products.
\end{ntt}

\subparagraph{Acknowledgements}
The authors gratefully acknowledge the hospitality and support
of Bielefeld University.
We are sincerely grateful to N.~Karpenko and A.~Vishik 
for the comments concerning rational cycles and norm varieties.
The first and the second
authors are sincerely grateful to A.~Bak for hospitality, advice, and the 
opportunity for all of us to work together.

\ 

{\small\sc
Steklov Mathematical Institute, St.~Petersburg, Russia

\ 

Fakult\"at f\"ur Mathematik, Universit\"at Bielefeld, Germany

\ 

Max-Planck-Institut f\"ur Mathematik, Bonn, Germany

\hspace{3em} {\tt kirill@math.uni-bielefeld.de}
}

\newpage

\bibliographystyle{chicago}

\end{document}